\documentclass[a4paper,10pt]{article}
\usepackage{amsmath}
\usepackage{amssymb}
\usepackage{graphics}
\usepackage{rotating}
\usepackage{cite}
\usepackage{amsfonts}
\usepackage{color}

\newcommand{\field}[1]{\mathbb{#1}}
\newtheorem{definition}{Definition}
\newtheorem{lemma}{Lemma}
\newtheorem{theorem}{Theorem}
\newtheorem{proposition}{Proposition}
\newtheorem{property}{Property}
\newtheorem{corollary}{Corollary}
\title{The quaternionic Cullen-regular product for a larger class of functions}
\author{Daniel Alay\'{o}n-Solarz}

\begin{document}

\maketitle

\begin{abstract}
We introduce the regular product for Cullen-regular quaternionic functions in a manner that does 
not depend upon a representation in power series but upon another, weaker kind of representation. The special case when the functions are represented as quaternionic power series is studied. We show that the regular ring of  quaternionic power series is a subring of the regular associative ring of real-analytical Hyperholomorphic functions.

\end{abstract}

\section{Introduction}

In \cite{GS2} G. Gentili and C. Stoppato have introduced the notion of the Cullen-regular product in two different ways; one way requires the pair of C-regular functions $f$ to admit a power series representation centered at the origin. Let $f$ and $g$ be a pair of quaternionic functions given by
\begin{displaymath}
f(p) : = \sum_{k = 0}^{\infty} p^n a_{n}, \ \ \ \  g(p): =  \sum_{k = 0}^{\infty} p^n b_{n}
\end{displaymath}
then the \textit{regular product}, $f \ast g (p)$  is given by
\begin{displaymath}
f \ast g (p) : = \sum_{n = 0}^{\infty} p^n c_{n}
\end{displaymath}
where
\begin{displaymath}
c_{n} :=  \sum_{k = 0}^{n} a_{k} b_{n-k}.
\end{displaymath}
Since functions defined as power series are Cullen-regular we obtain a new Cullen-regular function. The problem with this first kind of definition is that not all Cullen-regular functions can be expressed as such series, so this regular product is not well-defined for all Cullen-regular functions. \\
The second approach for the definition of the regular product is given by the following closed formula:
 \begin{displaymath}
 f \ast g(p):= f(p)g(f(p)^{-1} \ p \ f(p)),
 \end{displaymath}
 which has the advantage of being defined for every pair of functions $f$ and $g$. In \cite{GS2} it is proved that for power series this remarkable simple closed formula coincides with the regular product. We consider likely that this closed formula is in fact the regular product but, to our best knowledge, this has not been proved in general.
\\
In this paper we will introduce and study a third, independent approach to the problem of how to endow the Cullen-regular functions with a Cullen-regular product.
The main advantage of this approach is that is based in a weaker type of representation than power series. It turns out this representation must exist if the function is of class $C^{2}$. The problem with our definition is that associativity does not immediately follows and it is defined in terms of a first-order differential operator, so the result is a function of lesser differentiability class. We tackle this problem by restricting this product to real-analytic hyperholomorphic functions. In this case the product turns out to be associative. We discuss why a hyperholomorphic function of class $C^3$ is real-analytical. Then an inverse for this product is shown to exist if the function is real-analytical hyperholomorphic. 
\section{Preliminaries}
Let $\field{H}$ denote the quaternions and let $p$ be a quaternion. Thus $p$ can be written as:
\begin{displaymath}
p: = t + xi + yj + zk.
\end{displaymath}
Let $\iota$ denote the following function:
\begin{displaymath}
\iota := \frac{xi + yj + zk}{\sqrt{x^2+y^2+z^2}}
\end{displaymath}
Note that for all quaternion $p$ we have $\iota^2 = -1$. Now setting $r: =\sqrt{x^2+y^2+z^2}$. A quaternion  $p$ can be written as:
\begin{displaymath}
p = t + r \iota.
\end{displaymath}
A quaternionic valued function $f$ of one quaternionic variable $p$ of class $C^1$ is called \textbf{Cullen-regular} or C-regular if the expression
\begin{displaymath}
(\frac{\partial}{\partial t} + \iota \frac{\partial}{\partial r})f 
\end{displaymath}
vanishes identically.

\section{Complex-like regular functions}
 
We start restricting our attention to a special class of C-regular functions.
\begin{definition}
A quaternionic function  $f$ is said to admit a \textbf{complex-like regular form} if there exists quaternionic functions $u$ and $v$ of class $C^1$ such that
\begin{itemize}
\item $f = u + \iota v$
\item $u$ and $v$ satisfy the Cauchy-Riemann equations
\begin{displaymath}
\frac{\partial u}{\partial t}=\frac{\partial v}{\partial r},
\end{displaymath}
\begin{displaymath}
\frac{\partial u}{\partial r}=-\frac{\partial v}{\partial t}.
\end{displaymath} 
\end{itemize}
\end{definition}
What we have defined is a stronger form of C-regularity:
\begin{proposition}
Suppose a function $f$ admits a complex-like regular form. Then $f$ is C-regular.
\end{proposition}
Given two functions admitting a complex-like regular form we can always define a third, C-regular function:
\begin{proposition} Suppose $f_1$ and $f_2$ admit complex-like regular forms such that $f_1 = u_1 + \iota v_1$, $f_2 = u_2 + \iota v_2$ then the function defined as:
\begin{displaymath}
u_{1}u_{2} -v_{1}v_{2} +  \iota(v_{1}u_{2} +u_{1}v_{2}),
\end{displaymath}
is C-regular.
\end{proposition}
\textbf{Proof} It suffices to show that the complex-like function defined satisfies the Cauchy-Riemann equations:
\begin{displaymath}
\frac{\partial }{\partial t} (u_{1}u_{2} -v_{1}v_{2})=\frac{\partial }{\partial r}(v_{1}u_{2} +u_{1}v_{2}),
\end{displaymath}
\begin{displaymath}
\frac{\partial }{\partial r}(u_{1}u_{2} -v_{1}v_{2})=-\frac{\partial }{\partial t}(v_{1}u_{2} +u_{1}v_{2}).
\end{displaymath} 
which is a rather simple calculation. \\
\begin{proposition}
Suppose the function $f$ admits a complex-like regular form $f = u + \iota v$, where $u$ and $v$ are real valued. Then the quaternionic inverse of $f$ is C-regular.
\end{proposition}
\textbf{Proof} Let $f= u + \iota v$ be such that $u$ and $v$ are real-valued and satisfy the Cauchy-Riemann equations. Then its quaternionic inverse is:
\begin{displaymath}
f^{-1} = \frac{u}{u^2+v^2} + \iota (\frac{-v}{u^2+v^2} )
\end{displaymath}
it suffices to show then that if $u$ and $v$ satisfy the Cauchy-Riemann equations as vectors then:
\begin{displaymath}
\frac{\partial}{\partial t}( \frac{u}{u^2+v^2})=\frac{\partial}{\partial r}(\frac{-v}{u^2+v^2} )
\end{displaymath}
\begin{displaymath}
\frac{\partial}{\partial r}(\frac{u}{u^2+v^2})=-\frac{\partial}{\partial t}(\frac{-v}{u^2+v^2} )
\end{displaymath} 
which is again a rather simple calculation.\\
Observe that the complex-like regular form is not unique. So we don't have a well-defined regular product yet. In principle there is no unicity for a complex-like regular form. Consider the function $f=1$, which can be written as: $1 +\iota(0)$ or $0+ \iota(- \iota) $. However, if one has a systematic manner to write arbitrary C-regular in a complex-like regular form then one would also be able to systematically construct a rule asigning to a pair of C-regular functions a third C-regular function.    
\section{A characterization of $C^{2}$ C-regular functions}
\begin{theorem}
A quaternionic function of class $C^{2}$ $f$ is C-regular if and only $f$ admits a complex-like regular form.
\end{theorem}
\textbf{Proof} If there exists $u$ and $v$ quaternionic functions satisfying the Cauchy-Riemann equations then  $f := u + \iota v$ is C-regular. If $f$ is C-regular then:
\begin{displaymath}
\frac{\partial f}{\partial t}  + \iota \frac{\partial f}{\partial r} = 0
\end{displaymath}
Now parametrizing $\iota$ as
\begin{displaymath}
\iota = (\cos\alpha \sin\beta, \sin\alpha \sin\beta, \cos \beta).
\end{displaymath}
Now we define the following operator:
\begin{equation*}
\frac{\partial}{\partial \iota} := ({\iota}_{\alpha})^{-1}\frac{\partial}{\partial \alpha} +  ({\iota}_{\beta})^{-1}\frac{\partial}{\partial \beta}.
\end{equation*}
Which has the following fundamental property:
\begin{displaymath}
\frac{\partial}{\partial \iota} (\iota f) + \iota \frac{\partial}{\partial \iota} (f) +  = 2f
\end{displaymath}
Now we have:
\begin{displaymath}
\frac{\partial}{\partial \iota} (\frac{\partial f}{\partial t}  + \iota \frac{\partial f}{\partial r}) = 0
\end{displaymath}
using the fundamental property of the operator $\frac{\partial}{\partial \iota}$:
\begin{displaymath}
\frac{\partial}{\partial \iota} \frac{\partial f}{\partial t} + 2  \frac{\partial f}{\partial r} - \iota \frac{\partial}{\partial \iota} \frac{\partial f}{\partial r} = 0
\end{displaymath}
note that this expression exists because $f$ is of class $C^{2}$. Now setting:
\begin{displaymath}
u := \frac{1}{2} \frac{\partial}{\partial \iota}(\iota f),
\end{displaymath}
\begin{displaymath}
v := \frac{1}{2} \frac{\partial}{\partial \iota}(f).
\end{displaymath}
We see then that:
\begin{displaymath}
\frac{\partial v}{\partial t} +   \frac{\partial f}{\partial r} -\iota \frac{\partial v}{\partial r} = 0
\end{displaymath}
but this occurs if and only if:
\begin{displaymath}
\frac{\partial v}{\partial t} + \frac{\partial u}{\partial r} = 0,
\end{displaymath}
we then repeat this procedure for the expression:
\begin{displaymath}
\iota \frac{\partial f}{\partial t}  -  \frac{\partial f}{\partial r} = 0
\end{displaymath}
in order to obtain:
\begin{displaymath}
\frac{\partial u}{\partial t}  - \frac{\partial v}{\partial r}= 0.
\end{displaymath}
Because the complex-like regular form appearing in the characterization Theorem is important for our needs we shall distinguish it:
\begin{definition}
A  function $f$ is said to admit a \textbf{proper complex-like regular form} if the functions defined as:
\begin{displaymath}
u := \frac{1}{2} \frac{\partial}{\partial \iota} (\iota f)
\end{displaymath}
and
\begin{displaymath}
v := \frac{1}{2} \frac{\partial}{\partial \iota} ( f)
\end{displaymath}
satisfy the Cauchy-Riemann equations:
\begin{displaymath}
\frac{\partial u}{\partial t}=\frac{\partial v}{\partial r},
\end{displaymath}
\begin{displaymath}
\frac{\partial u}{\partial r}=-\frac{\partial v}{\partial t}.
\end{displaymath} 
\end{definition}
We can re-state the characterization theorem as:
\begin{corollary}
A $C^2$ quaternionic is C-regular if and only if admits a proper complex-like regular form. 
\end{corollary}
\section{The regular product}
\begin{definition}
Let $f_1$ and $f_2$ be a pair of $C^2$ C-regular functions, then their \textbf{regular product} $f \ast g$ is given by:
\begin{displaymath}
f_1 \ast f_2 := (u_1 u_2 - v_1 v_2) + \iota(u_1 v_2 + v_1 u_2)
\end{displaymath} 
where
\begin{displaymath}
u_1 := \frac{1}{2} \frac{\partial}{ \partial \iota}(\iota f_1),
\end{displaymath}
\begin{displaymath}
v_1 := \frac{1}{2} \frac{\partial}{ \partial \iota}(f_1),
\end{displaymath}
\begin{displaymath}
u_2 := \frac{1}{2} \frac{\partial}{ \partial \iota}(\iota f_2),
\end{displaymath}
\begin{displaymath}
v_2 := \frac{1}{2} \frac{\partial}{ \partial \iota}(f_2)
\end{displaymath}
\end{definition}
Notice how we fix the regular-product for the proper complex-like regular form appearing in the characterization. \\
\textbf{Remark}: Some information is lost in the regular product: the complex-like regular form is obtained from two \textbf{proper} complex-like regular forms. However the quality of  'being proper' might not be present in the form givem by the product. Further hipothesis are needed. 
\begin{property}
Let $f$ and $g$ be $C^{2}$ and C-regular. Then their regular product satisfy:
\begin{itemize}
\item $f \ast g$ is C-regular. 
\item $ f \ast (g + h)$ = $ f \ast g + f \ast h$ and $(f + g) \ast h = f \ast h + g \ast h$ (\textbf{Distribuitivity})
\item $f \ast 1 =  f \ast 1 = f$ (\textbf{Neutral element})
\end{itemize}
\end{property}
\begin{property}
Suppose the C-regular function $f$ admits a proper complex-like regular form $f = u+ \iota v$ such that $u$ and $v$ are real valued. Then, for every $C^2$ C-regular function $g$ we have
\begin{displaymath}
f \ast g = g \ast f = fg.
\end{displaymath}
\end{property}
Which means that the regular product and the quaternionic product can coincide. Following Gentili and Stoppato we define the \textbf{regular conjugate} denoted by $f^{c}$ in the following way
\begin{displaymath}
f^{c} := \bar{u} + \iota \bar{v}
\end{displaymath}
and its \textbf{symmetrization}
\begin{displaymath}
f^{s} :=  |u|^2 - |v|^2 + \iota (u \bar v + v \bar u) =  |u|^2 - |v|^2 + \iota (\bar u v + \bar v u)
\end{displaymath}
Note the following identity:
\begin{displaymath}
\bar u v + \bar v u = u \bar v + v \bar u = 2 \ (u \cdot v)
\end{displaymath}
where the right hand-side is the common euclidian inner-product of $u$ and $v$ as $4$-vectors.
\begin{proposition}
The symmetrization $f^{s}$ of a C-regular function is a C-regular function that admits a complex-like form $f^{s}=u + \iota v$ such that $u$ and $v$ are real-valued functions.
\end{proposition}
\textbf{Proof} Take $f_1 = f$ and $f_2 :=f^c$ in \textbf{Proposition 2} 
\begin{displaymath}
u_{1}u_{2} - v_{1}v_{2} + \iota(u_{1}v_{2} + v_{1}u_{1}) = u \bar u - v \bar v + \iota(u \bar v + v \bar u)
\end{displaymath}
Because $u$ and $v$ satisfy the Cauchy-Riemann equations as vectors it follows that the pairs $\bar{u}$ and $\bar{v}$ satisfy this equations also. This shows that $f^{c}$ admits a complex-like regular form. Because the regular product of two functions admitting a complex-like regular form is a function that admits a complex-like regular form it follows that $f^{s}$ admits a complex-like regular form. 
As a consequence:
\begin{proposition}
The function $\frac{1}{f^{s}}$, where defined, admits a complex-like regular form with real-valued functions $u$ and $v$.
\end{proposition}
So now we can define, for this special case the \textbf{regular reciprocal} or reciprocal, of a C-regular function that admits a complex-like representation.
\begin{displaymath}
f^{- \ast} : = \frac{1}{f^{s}} \ast f^{c} = f^{c} \ast \frac{1}{f^{s}} = \frac{1}{f^{s}} f^{c}  .
\end{displaymath}
In order to ensure that the reciprocal is the inverse under the regular product one needs associativity. We remark that some information is lost with the regular product. 
 We will deal with this later on.
\section{Special case: The regular product for quaternionic power series}

Let $a$ be a constant. Observe that in this case the characterization of \textbf{Theorem 1} gives the following:
\begin{displaymath}
u = a
\end{displaymath}
and
\begin{displaymath}
v = 0.
\end{displaymath}
This implies that for constants $a$ and $b$ the regular product coincides with the quaternionic product:
\begin{displaymath}
a \ast b = ab.
\end{displaymath}
Now the case of the power function $p^n$ for integers $n$, inductively one shows that $p^{n}$ can be written as $u_{n} + \iota v_{n}$ where $u_{n}$ and $v_{n}$ are real-valued functions that only depend on the variables $t$ and $r$, and therefore, as $C^{2}$ C-regular functions the expression $u_{n} + \iota v_{n}$ coincides with the $u$ and $v$ given in \textbf{Theorem 1}. Therefore, is straight-forward to check that
\begin{displaymath}
p^{n} \ast p^{m} = p^{m} \ast p^{n} = p^{n+m},
\end{displaymath}
and
\begin{displaymath}
a \ast p^{n} = p^{n} \ast a = p^{n}a.
\end{displaymath}
In general we can show that:
\begin{displaymath}
(p^n a) \ast (p^m b) = p^{n+m}ab
\end{displaymath}
For this, setting $p^{n} a$ in proper complex-like regular form we have $p^n a = u_{n} a + \iota v_{n} a $ and $p^{m} = u_{m} b + \iota v_{m} b$. Then
\begin{displaymath}
(p^n a) \ast (p^m b) = u_{n} a u_{m} b - v_{n} a v_{m} b + \iota (v_{n} a u_{m} b + u_{n} a v_{m} b)
\end{displaymath}
Now the key part is to observe that the functions $u_{n,m}$, $v_{n,m}$ are all real valued, so the right hand side can be rewritten as, since the quaternionic product is associative:
\begin{displaymath}
\Big ( u_{n}  u_{m}  - v_{n}  v_{m}  + \iota (v_{n}  u_{m}  + u_{n}  v_{m} ) \Big ) ab = p^{n+m}ab.
\end{displaymath}

in view of this, let $f$ and $g$ be two functions constructed as power series,
\begin{displaymath}
f(p) : = \sum_{n = 0}^{\infty} p^n a_{n}, \ \ \ \  g(p): =  \sum_{n = 0}^{\infty} p^n b_{n}
\end{displaymath}
then
\begin{displaymath}
f \ast g = (\sum_{n = 0}^{\infty} p^n a_{n}) \ast (\sum_{n = 0}^{\infty} p^n b_{n})
\end{displaymath}
using, linearity is made of expressions of the form $(p^{k} a_{k}) \ast (p^{n-k}  b_{n-k})$, since all of this expressions can be reduced to $p^{n}a_{n-k}b_{k}$ we have that the regular product must be:
\begin{displaymath}
f \ast g  : = \sum_{n = 0}^{\infty} p^n c_{n}
\end{displaymath}
where
\begin{displaymath}
c_{n} :=  \sum_{k = 0}^{n} a_{k} b_{n-k}.
\end{displaymath}
For a quaternonic constant $a$ the regular conjugate coincides with the quaternionic conjugate, so
\begin{displaymath}
a^{c} = \bar {a}
\end{displaymath}
But it leaves the function $p^{n}$ invariant:
\begin{displaymath}
(p^n)^c = p^n
\end{displaymath}
We also have
\begin{displaymath}
(p^n a)^{c} = p^n \bar{a} 
\end{displaymath}
Since $(f + g)^c = f^c + g^c$ we must conclude that for a function $f$ defined as a power series
\begin{displaymath}
f(p) : = \sum_{n = 0}^{\infty} p^n a_{n}
\end{displaymath}
We have that its regular conjugate is:
\begin{displaymath}
f^c =  \sum_{n = 0}^{\infty} p^n \bar{a_{n}}.
\end{displaymath}
In other words, in the case of the power series our definition of the regular product coincides with the definition given by Gentili and Stoppato. Functions constructed as such quaternionic power series are examples of Hyperholomorphic functions.

\section{Real-analytical Hyperholomorphic functions as an Associative regular ring}
A regular ring is a sub-class of the C-regular functions such that is closed under the regular product. One example is the case of the power series. Gentili and Stoppato have shown, see \cite{GS2} that for this sub-ring the
regular product turns out to be associative. Associative subrings have the advantage that the regular reciprocal can be defined and is the inverse under the regular product:
\begin{proposition} Let $f$ be a function defined as a power series, then
\begin{displaymath}
f^{-\ast} \ast f = f \ast f^{-\ast} = 1
\end{displaymath}
\end{proposition}  
\textbf{Proof}. Since this subring is associative, see \cite{GS2} we have the following equalities:
\begin{displaymath}
f^{-\ast} \ast f = (\frac{1}{ f^{s}} \ast f^{c}) \ast f = \frac{1}{ f^s} \ast (f^{c} \ast f) = \frac{1}{f^s}  \ast f^{s} =  \frac{1}{f^s} f^{s} =1.
\end{displaymath}
and
\begin{displaymath}
f \ast f^{-\ast} =   f \ast (\frac{1}{ f^{s}} \ast f^{c})  = f \ast ({ f^{c}} \ast \frac{1}{f^{s}}) = (f \ast f^ c) \ast \frac{1}{f^s} = f^s  \frac{1}{f^s} =1.
\end{displaymath}
We will now show for every triad of hyperholomorphic functions $f$, $g$ and $h$ hyperholomorphic quaternionic we have
\begin{displaymath}
f \ast (g \ast h) = (f \ast g) \ast h
\end{displaymath} 
First we recall the definition of a hyperholomorphic function:
\begin{definition}
A C-regular function $f$ of class $C^2$ is said to be hyperholomorphic if its proper complex-like regular form $u + \iota v$ satisfy the following
modified Cauchy-Riemann system:
\begin{eqnarray}
\frac{\partial v}{\partial \alpha}(\sin\beta)^{-1}+\frac{\partial u}{\partial \beta} = 0,
\\
\frac{\partial u}{\partial \alpha}(\sin\beta)^{-1}-\frac{\partial v}{\partial \beta}= 0.
\end{eqnarray}
\end{definition}
We want to prove the following important property of the Hyperholomorphic functions:
\begin{property}
The regular product of two hyperholomorphic functions is in a proper regular complex-like form.
\end{property}
For this we want to characterize proper complex-like regular forms:
\begin{lemma}
Let $f$ be a $C^2$ C-regular function. Let $u + \iota v$ be the proper complex-like form of $f$. Then $u$ and $v$ satisfy the following equation:
\begin{displaymath}
\frac{\partial u }{\partial \iota} = \iota \frac{\partial v }{\partial \iota}
\end{displaymath}
\end{lemma}
\textbf{Proof} Observe that, using the definition of $v$ recursively:
\begin{displaymath}
2v = \frac{\partial f}{\partial \iota} = \frac{\partial}{ \partial \iota}(u + \iota v) =  \frac{\partial u}{ \partial \iota} +2v - \iota  \frac{\partial v}{ \partial \iota}
\end{displaymath}
Which is a compatibility condition. A rather lenghty calculation is needed to the following
\begin{proposition}
$f_1 = u_1 + \iota v_1$ and $f_2 = u_2 + \iota v_2$ are hyperholomorphic then:
\begin{displaymath}
\frac{\partial}{\partial \iota}(u_{1}u_{2} -v_{1}v_{2}) = \iota \frac{\partial}{\partial \iota}(v_{1}u_{2} +u_{1}v_{2})
\end{displaymath}
\end{proposition}
Now assuming this proposition holds we can show that
\begin{proposition}
Let $f_1$ and $f_2$ be hyperholomorphic then their regular product satisfy the following equality:
\begin{displaymath}
\frac{\partial}{\partial \iota}(f_1 \ast f_2 ) =   \frac{1}{2} \Big ( \frac{\partial}{\partial \iota}(\iota f_1 )\frac{\partial}{\partial \iota}(f_2 )  + \frac{\partial}{\partial \iota}(f_1) \frac{\partial}{\partial \iota}(\iota f_2 ) \Big )
\end{displaymath}
\end{proposition}
\textbf{Proof} The proof is the same as the \textbf{Lemma 1}:
\begin{displaymath}
\frac{\partial}{\partial \iota}(f_1 \ast f_2 ) =  \frac{\partial}{\partial \iota}(u_{1}u_{2} -v_{1}v_{2}+ \iota(v_{1}u_{2} +u_{1}v_{2}) )
\end{displaymath}
\begin{displaymath}
= \frac{\partial}{\partial \iota}(u_{1}u_{2} -v_{1}v_{2}) + 2 (v_{1}u_{2} +u_{1}v_{2}) -\iota  \frac{\partial}{\partial \iota}(v_{1}u_{2} +u_{1}v_{2})
\end{displaymath}
We also have that the regular product is compatible with the modified Cauchy-Riemann equations that rule hyperholomorphic functions:
\begin{proposition}
Let $f_1$ and $f_2$ be hyperholomorphic then:
\begin{eqnarray*}
(\sin\beta)^{-1} \frac{\partial }{\partial \alpha} (v_{1}u_{2} +u_{1}v_{2} )+\frac{\partial }{\partial \beta} (u_{1}u_{2} -v_{1}v_{2}) = 0,
\\
(\sin\beta)^{-1} \frac{\partial }{\partial \alpha} (u_{1}u_{2} -v_{1}v_{2}) -\frac{\partial }{\partial \beta} (v_{1}u_{2} +u_{1}v_{2}) = 0.
\end{eqnarray*}
where
\begin{displaymath}
u_1 := \frac{1}{2} \frac{\partial}{ \partial \iota}(\iota f_1),
\end{displaymath}
\begin{displaymath}
v_1 := \frac{1}{2} \frac{\partial}{ \partial \iota}(f_1),
\end{displaymath}
\begin{displaymath}
u_2 := \frac{1}{2} \frac{\partial}{ \partial \iota}(\iota f_2),
\end{displaymath}
\begin{displaymath}
v_2 := \frac{1}{2} \frac{\partial}{ \partial \iota}(f_2)
\end{displaymath}
\end{proposition}
We can now state the following result:
\begin{proposition}
Let $f_1$ and $f_2$ be hyperholomorphic functions of class $C^3$. Then their regular product is an hyperholomorphic function.
\end{proposition}
as a consequence:
\begin{proposition}
Let $f_1$, $f_2$ and $f_3$ be real-analytical hyperholomorphic functions. Then
\begin{displaymath}
 (f_1 \ast f_2) \ast f_3 = f_1 \ast (f_2 \ast f_3)
\end{displaymath}
\end{proposition}
First we have:
\begin{displaymath}
f_1 \ast f_2 = u_{1}u_{2} -v_{1}v_{2}+ \iota(v_{1}u_{2} +u_{1}v_{2})
\end{displaymath}
Now to calculate $(f_1 \ast f_2) \ast f_3$ using the definition of the regular product we need to calculate:
\begin{displaymath}
\frac{\partial}{\partial \iota}((f_1 \ast f_2))
\end{displaymath}
and
\begin{displaymath}
\frac{\partial}{\partial \iota}(\iota (f_1 \ast f_2))
\end{displaymath}
But this equals:
\begin{displaymath}
\frac{\partial}{\partial \iota}(f_1 \ast f_2) =2( v_{1}u_{2} +u_{1}v_{2})
\end{displaymath}
and
\begin{displaymath}
\frac{\partial}{\partial \iota}(\iota (f_1 \ast f_2)) = 2(u_{1}u_{2} -v_{1}v_{2})
\end{displaymath}
Therefore
\begin{displaymath}
(f_1 \ast f_2 ) \ast f_3 = (u_{1}u_{2} -v_{1}v_{2})u_{3} -  (v_{1}u_{2} +u_{1}v_{2})v_{3} + \iota \Big ( (u_{1}u_{2} -v_{1}v_{2}) v_{3} +  ( v_{1}u_{2} +u_{1}v_{2}) u_{3}   \Big )
\end{displaymath}
but since the quaternionic product appearing in the right  hand side is associative we can get rid of the parenthesis to obtain:
\begin{displaymath}
 u_{1}u_{2}u_{3} -v_{1}v_{2}u_{3} -  v_{1}u_{2}v_{3} +u_{1}v_{2}v_{3} + \iota \Big ( u_{1}u_{2}v_{3} -v_{1}v_{2} v_{3} +   v_{1}u_{2}u_{3} +u_{1}v_{2} u_{3}   \Big )
 \end{displaymath}
which is equal to:
\begin{displaymath}
f_1 \ast( f_2  \ast f_3)
\end{displaymath}
Given $f$ a $C^2$ hyperholomorphic function in proper regular complex-like form $ f = u + \iota v$.
Recall the regular conjugate is $f^c := \bar u + \iota \bar v$.
\begin{proposition}
The regular conjugate of a $C^3$ hyperholomorphic function is an hyperholomorphic function. 
\end{proposition}
\textbf{Proof} First we want to show that the compatibility condition is met:
\begin{displaymath}
\frac{\partial \bar u}{\partial \iota} = \iota \frac{\partial \bar v}{\partial \iota},
\end{displaymath}
first note the that the pair $\bar u, \bar v$ still satisfy the same modified Cauchy-Riemann:
\begin{eqnarray}
\frac{\partial \bar v}{\partial \alpha}(\sin\beta)^{-1}+\frac{\partial \bar u}{\partial \beta} = 0,
\\
\frac{\partial\bar u}{\partial \alpha}(\sin\beta)^{-1}-\frac{\partial \bar v}{\partial \beta}= 0.
\end{eqnarray}
applying the definition, then one replaces every partial derivative in $u$ by its counterpart in $v$
\begin{displaymath}
\frac{\partial \bar u}{\partial \iota} =  ({\iota}_{\alpha})^{-1}\frac{\partial \bar u}{\partial \alpha} +  ({\iota}_{\beta})^{-1}\frac{\partial \bar u}{\partial \beta} =  ({\iota}_{\alpha})^{-1}\sin \beta \frac{\partial \bar v}{\partial \beta} +  ({\iota}_{\beta})^{-1} \frac{-1}{\sin \beta} \frac{\partial  \bar v}{\partial \alpha} 
\end{displaymath}
and now use the following identities:
\begin{displaymath}
 \frac{-1}{\sin \beta} (\iota_{\beta})^{-1}
 = \iota (\iota_{\alpha})^{-1}
\end{displaymath}
and
\begin{displaymath}
 (\iota_{\alpha})^{-1} \sin \beta = \iota (\iota_{\beta})^{-1}
\end{displaymath}
we see that the right hand side is:
\begin{displaymath}
\iota \frac{\partial \bar v}{\partial \iota}.
\end{displaymath}
This means that the pair $\bar u$ and $\bar v$ define a proper complex-like regular form and that this regular form satisfy the hyperholomorphicity condition. This implies that the symmetrization of a hyperholomorphic,
function defined as:
\begin{displaymath}
f^s : = f^c \ast f
\end{displaymath}
is necessarily hyperholomophic.
Observe that:
\begin{theorem}
{$C^3$ Hyperholomorphic fuctions are real-analytical}.
\end{theorem}
First we show the following:
\begin{proposition}
Let $f$ be a $C^3$ hyperholomorphic function, then $f$ satisfies:
\begin{displaymath}
D_{l} \Delta f = 0.
\end{displaymath}
\end{proposition}
\textbf{Remark} The reader will recognize this is Fueter's Theorem. \\
\textbf{Proof}
The operators $D_{l}$ and $\Delta$ commute with each other. So it suffices to show that $\Delta D_{l} f = 0$. Now:
\begin{displaymath}
D_{l}f = - \frac{2v}{r}
\end{displaymath}
where 
\begin{displaymath}
v = \frac{1}{2} \frac{\partial f}{\partial \iota}
\end{displaymath}
Since $v$ appears in the Cauchy-Riemann equations as a vector. Then all components of $v$ must be a zero of
\begin{displaymath}
\frac{\partial^2}{\partial t^2} +  \frac{\partial^2}{\partial r^2}
\end{displaymath}
since $v$ also appears in the modified Cauchy-Riemann equations then all of its components must be a zero of the angular part of the Laplacian:
\begin{displaymath}
 \frac{1}{\sin^2 \beta }\frac{\partial^2 }{\partial \alpha^2} +\frac{\partial^2 }{\partial \beta^2} + \cot \beta \frac{\partial }{\partial \beta}
\end{displaymath}
since the latter can be factorized as:
\begin{displaymath}
 \frac{1}{\sin^2 \beta }\frac{\partial^2 }{\partial \alpha^2} +\frac{\partial^2 }{\partial \beta^2} + \cot \beta \frac{\partial }{\partial \beta} = \frac{1}{\sin^2 \beta}(\frac{\partial}{\partial \alpha} - i \sin \beta \frac{\partial}{\partial \beta})(\frac{\partial}{\partial \alpha} + i \sin \beta \frac{\partial}{\partial \beta})
\end{displaymath}
This means that for a real component $v_{n}$, $n = 0,1,2,3$ there exists a correspondent real function $u_{n}$ such that the modified Cauchy-Riemann equations are satisfied. But then the complex function $u_{n}+i v_{n}$ must be a zero of:
\begin{displaymath}
(\frac{\partial}{\partial \alpha} - i \sin \beta \frac{\partial}{\partial \beta})
\end{displaymath}
Combining this observations one can conclude that:
\begin{displaymath}
\Delta \frac{v}{r} = \Big (\frac{\partial^2}{\partial t^2} +  \frac{\partial^2}{\partial r^2} + \frac{2}{r}\frac{\partial}{\partial r} + \frac{1}{r^2}( \frac{1}{\sin^2 \beta }\frac{\partial^2 }{\partial \alpha^2} +\frac{\partial^2 }{\partial \beta^2} + \cot \beta \frac{\partial }{\partial \beta}) \Big )(\frac{v}{r}) = 0.
\end{displaymath}
Therefore all components of the function $\frac{v}{r}$ are harmonical. Every harmonical real function is the real part of some Fueter-regular function and thus is real-analytical. Therefore the function $\frac{v}{r}$ is real-analytical. Since the function $r$ is also real-analytical (except at the real numbers in the quaternions)  then the product of them, namely, $v$ is real-analytical. We can repeat this result to the function $\iota f$ in order to obtain that $u$ is also real-analytical. Since $f$ is the sum of the real-analytical functions $u$ and $\iota v$ then we conclude $f$ is also real-analytical except possibly at the real numbers. If the C-regular function is defined on some open subset of the real numbers then we can express $f$ in a small neighbourhood of a real numbers as a quaternionic power series, then $f$ is necessarily real-analytical. \\
Now suppose $f$ is a real-analytical hyperholomorphic function, then $f^c$ is also a real-analytical hyperholomorphic. Therefore the function $f^s := f \ast f^c = f^c \ast f$ must be analytical hyperholomorphic. Now denoting
\begin{displaymath}
\tilde{u} :=  |u|^2 - |v|^2
\end{displaymath} 
and
\begin{displaymath}
\tilde{v} := (u \bar v + v \bar u)
\end{displaymath}
so that $f = \tilde{u} + \iota \tilde{v}$. Note that $\tilde{u}$ and $\tilde{v}$ are real valued, and satisfy:
\begin{displaymath}
\frac{\partial \tilde{u}}{\partial t}=\frac{\partial \tilde{v}}{\partial r},
\end{displaymath}
\begin{displaymath}
\frac{\partial \tilde{u}}{\partial r}=-\frac{\partial \tilde{v}}{\partial t}.
\end{displaymath} 
and
\begin{eqnarray}
\frac{\partial \tilde{v}}{\partial \alpha}(\sin\beta)^{-1}+\frac{\partial \tilde{u}}{\partial \beta} = 0,
\\
\frac{\partial \tilde{u}}{\partial \alpha}(\sin\beta)^{-1}-\frac{\partial \tilde{v}}{\partial \beta}= 0.
\end{eqnarray}
Because $\tilde{u}$ and $\tilde{v}$ are real-valued then the quaternionic inverse of this function, is
\begin{displaymath}
\frac{\tilde{u}}{\tilde{u}^2+\tilde{v}^2} + \iota (\frac{-\tilde{v}}{\tilde{u}^2+\tilde{v}^2} )
\end{displaymath}
which we now rewrite as $\hat{u}$ and $\hat{v}$. We have again this pair satisfy the Cauchy-Riemann and modified Cauchy-Riemann. From this we can show then that the pair must does satisfy the compatibility condition:
\begin{displaymath}
\frac{\partial \hat{u}}{\partial \iota} = \iota \frac{\partial \hat{v}}{\partial \iota}. 
\end{displaymath}
Therefore if $\hat{u}$,$\hat{v}$
are of class $C^2$, the function $\hat{u} + \iota \hat{v}$ is hyperholomorphic. Now setting
\begin{displaymath}
f^{- \ast} : = \frac{1}{f^{s}} \ast f^{c} = f^{c} \ast \frac{1}{f^{s}} = \frac{1}{f^{s}} f^{c}.
\end{displaymath}
Because the right right hand-side are hyperholomorphics functions, then the regular reciprocal $f^{- \ast}$, if defined, is a hyperholomorphic function. Notice that we have assumed the symmetrization $f^{s}$ to be
non-zero in some suitable open set.  We can now state this result as:
\begin{theorem}
The class of real-analytic hyperholomorphic functions on a domain $\Omega$ is a regular, non-commutative but associative ring extension to the regular ring of power series.  If the symmetrization
$f^{s}$ is non-zero then  the hyperholomorphic function $f$ is a ring unit.
\end{theorem}





\end{document}